\documentclass{amsart}%
\usepackage{amsfonts}
\usepackage{amssymb,latexsym}
\usepackage{amsmath}
\usepackage{amssymb}
\newtheorem{theorem}{Theorem}

\newtheorem{lemma}[theorem]{Lemma}
\newtheorem{claim}[theorem]{Claim}
\theoremstyle{remark}
\newtheorem{remark}[theorem]{Remark}
\begin{document}
\title{Minimal clones with many majority operations}
\author{Mike Behrisch}
\address[Mike Behrisch]{Institut f\"{u}r Algebra, TU Dresden, D-01062 Dresden, Germany}
\email{mike.behrisch@mailbox.tu-dresden.de}
\author{Tam\'{a}s Waldhauser}
\address[Tam\'{a}s Waldhauser]{University of Luxembourg, Mathematics Research Unit, 6
rue Richard Coudenhove-Kalergi, L-1359 Luxembourg, Luxembourg, and\\
Bolyai Institute, University of Szeged, Aradi v\'{e}rtan\'{u}k tere 1, H-6720
Szeged, Hungary}
\email{twaldha@math.u-szeged.hu}

\begin{abstract}
We present two minimal clones containing 26 and 78 majority operations
respectively, more than any other previously known example.

\end{abstract}
\maketitle

\section{Introduction}

A \emph{clone} is a family of finitary operations defined on a set $A$ that is
closed under composition of functions and contains all projections (which will
also be called \emph{trivial functions}). Given a set $F$ of operations on
$A$, the functions obtained from elements of $F$ and from projections by means
of compositions form the smallest clone containing $F$. This is \emph{the
clone generated by }$F$, and we denote this clone by~$\left[  F\right]  $.
This clone is nothing else but the clone of term functions of the 
algebra~$\left(  A;F\right)  $.

The set of all clones on a given base set $A$ is a complete lattice; the
largest element of this lattice is the clone of all operations on $A$, and the
smallest element is the clone containing projections only. The latter is
called the \emph{trivial clone}, denoted by~$\mathcal{I}$. A \emph{minimal
clone} is an atom in the clone lattice, i.e., a nontrivial clone, whose only
proper subclone is $\mathcal{I}$. As opposed to the case of maximal clones
(coatoms of the clone lattice), the description of minimal clones is still an
open problem, although there are numerous partial results. Here we review only
those facts about minimal clones that we need in the sequel, for an overview
of minimal clones we refer the reader to the survey papers \cite{Csminicourse}
and \cite{Qsurv}; for general reference on clones see~\cite{Lau,PK,SzAclUA}.

It follows from the definition that every minimal clone is generated by any
one of its nontrivial members, and a nontrivial function $f$ generates a
minimal clone iff%
\begin{equation}
f\in\left[  h\right]  \text{ holds for all }h\in\left[  f\right]
\setminus\mathcal{I}. \label{minimality criterion}%
\end{equation}
We will consider clones generated by a \emph{majority operation}, i.e., by a
ternary operation $f$ satisfying%
\[
f\left(  a,a,b\right)  =f\left(  a,b,a\right)  =f\left(  b,a,a\right)
=a\text{ for all }a,b\in A.
\]
As it was shown in \cite{Cscons}, in this case all ternary functions in
$\left[  f\right]  $ are majority operations, except for the three projections.

For any clone $\mathcal{C}$, let $\mathcal{C}^{\left(  3\right)  }$ denote the
set of ternary operations belonging to $\mathcal{C}$. The composition of functions
yields a quaternary operation on $\mathcal{C}^{\left(  3\right)  }$ as we
have one outer function and three inner functions in a composition.
Furthermore, we may regard the three ternary projections in 
$\mathcal{C}^{\left(  3\right)  }$ as nullary operations. With these operations 
$\mathcal{C}^{\left(  3\right)  }$ becomes
an algebra of type $\left(  4,0,0,0\right)  $, a so-called \emph{unitary
Menger algebra} of rank $3$ (cf. \cite{Me}). Among the many pleasant
properties of majority operations, there is one that is especially useful in
the investigation of minimal clones: if $f$ is a majority operation, then the
minimality of the clone $\mathcal{C}=\left[  f\right]  $ is determined by its
ternary part, i.e., it suffices to check the minimality criterion
(\ref{minimality criterion}) only for \emph{ternary} functions (see
\cite{Cscons,Wafew}). Formally, if $f$ is a majority operation, then $\left[
f\right]  $ is a minimal clone iff%
\begin{equation}
f\in\left[  h\right]  \text{ holds for all }h\in\left[  f\right]  ^{\left(
3\right)  }\setminus\mathcal{I}. \label{minimality criterion maj}%
\end{equation}
If the base set is finite, this means that there are only finitely many
functions $h$ to be tested, hence, at least in principle, it can be done by computer.

There are very few examples of minimal clones generated by a majority
operation (while there is an abundance of examples of other types of minimal
clones). Two general examples are:\ the median function $\left(  x\wedge
y\right)  \vee\left(  y\wedge z\right)  \vee\left(  z\wedge x\right)  $ on any
lattice (see, e.g., \cite{PK}), and the dual discriminator function on any set
(see \cite{CsG,FP}). All other examples came from systematic investigations of
minimal clones on small sets. B.~Cs\'{a}k\'{a}ny determined all minimal clones
on any three-element set in \cite{Cs3all}, and among the clones he found there
are up to isomorphism three that are generated by a majority operation. These
clones contain 1, 3 and 8 majority operations, respectively; see
Table~\ref{table 3maj}. (Here, and in the other tables we omit those triples
where the majority rule determines the values of the functions.)
\begin{table}[t]
\centering%
\begin{tabular}
[c]{ccccccccccccc}\hline
\multicolumn{1}{|c}{} & \multicolumn{1}{|c}{$m_{1}$} &
\multicolumn{1}{|c}{$m_{2}$} &  &  & \multicolumn{1}{|c}{$m_{3}$} &  &  &  &
&  &  & \multicolumn{1}{c|}{}\\\hline
\multicolumn{1}{|c}{$(1,2,3)$} & \multicolumn{1}{|c}{$1$} &
\multicolumn{1}{|c}{$1$} & $2$ & $3$ & \multicolumn{1}{|c}{$3$} & $3$ & $1$ &
$3$ & $1$ & $1$ & $3$ & \multicolumn{1}{c|}{$1$}\\
\multicolumn{1}{|c}{$(2,3,1)$} & \multicolumn{1}{|c}{$1$} &
\multicolumn{1}{|c}{$2$} & $3$ & $1$ & \multicolumn{1}{|c}{$3$} & $1$ & $3$ &
$3$ & $1$ & $3$ & $1$ & \multicolumn{1}{c|}{$1$}\\
\multicolumn{1}{|c}{$(3,1,2)$} & \multicolumn{1}{|c}{$1$} &
\multicolumn{1}{|c}{$3$} & $1$ & $2$ & \multicolumn{1}{|c}{$3$} & $3$ & $3$ &
$1$ & $1$ & $1$ & $1$ & \multicolumn{1}{c|}{$3$}\\\hline
\multicolumn{1}{|c}{$(2,1,3)$} & \multicolumn{1}{|c}{$1$} &
\multicolumn{1}{|c}{$2$} & $1$ & $3$ & \multicolumn{1}{|c}{$1$} & $3$ & $1$ &
$1$ & $3$ & $1$ & $3$ & \multicolumn{1}{c|}{$3$}\\
\multicolumn{1}{|c}{$(1,3,2)$} & \multicolumn{1}{|c}{$1$} &
\multicolumn{1}{|c}{$1$} & $3$ & $2$ & \multicolumn{1}{|c}{$1$} & $1$ & $1$ &
$3$ & $3$ & $3$ & $3$ & \multicolumn{1}{c|}{$1$}\\
\multicolumn{1}{|c}{$(3,2,1)$} & \multicolumn{1}{|c}{$1$} &
\multicolumn{1}{|c}{$3$} & $2$ & $1$ & \multicolumn{1}{|c}{$1$} & $1$ & $3$ &
$1$ & $3$ & $3$ & $1$ & \multicolumn{1}{c|}{$3$}\\\hline
&  & $d_{1}$ & $d_{2}$ & $d_{3}$ &  &  &  &  &  &  &  & \\
&  &  &  &  &  &  &  &  &  &  &  &
\end{tabular}
\caption{{}Majority operations generating a minimal clone on the three-element
set $\left\{  1,2,3\right\}  $}%
\label{table 3maj}%
\end{table}\begin{table}[t]
\centering%
\begin{tabular}
[c]{ccccccccccccc}\hline
\multicolumn{1}{|c}{} & \multicolumn{1}{|c}{$M_{1}$} &
\multicolumn{1}{|c}{$M_{2}$} &  &  & \multicolumn{1}{|c}{$M_{3}$} &  &  &  &
&  &  & \multicolumn{1}{c|}{}\\\hline
\multicolumn{1}{|c}{$(1,2,3)$} & \multicolumn{1}{|c}{$4$} &
\multicolumn{1}{|c}{$4$} & $2$ & $3$ & \multicolumn{1}{|c}{$3$} & $3$ & $4$ &
$3$ & $4$ & $4$ & $3$ & \multicolumn{1}{c|}{$4$}\\
\multicolumn{1}{|c}{$(2,3,1)$} & \multicolumn{1}{|c}{$4$} &
\multicolumn{1}{|c}{$2$} & $3$ & $4$ & \multicolumn{1}{|c}{$3$} & $4$ & $3$ &
$3$ & $4$ & $3$ & $4$ & \multicolumn{1}{c|}{$4$}\\
\multicolumn{1}{|c}{$(3,1,2)$} & \multicolumn{1}{|c}{$4$} &
\multicolumn{1}{|c}{$3$} & $4$ & $2$ & \multicolumn{1}{|c}{$3$} & $3$ & $3$ &
$4$ & $4$ & $4$ & $4$ & \multicolumn{1}{c|}{$3$}\\\hline
\multicolumn{1}{|c}{$(2,1,3)$} & \multicolumn{1}{|c}{$4$} &
\multicolumn{1}{|c}{$2$} & $4$ & $3$ & \multicolumn{1}{|c}{$4$} & $3$ & $4$ &
$4$ & $3$ & $4$ & $3$ & \multicolumn{1}{c|}{$3$}\\
\multicolumn{1}{|c}{$(1,3,2)$} & \multicolumn{1}{|c}{$4$} &
\multicolumn{1}{|c}{$4$} & $3$ & $2$ & \multicolumn{1}{|c}{$4$} & $4$ & $4$ &
$3$ & $3$ & $3$ & $3$ & \multicolumn{1}{c|}{$4$}\\
\multicolumn{1}{|c}{$(3,2,1)$} & \multicolumn{1}{|c}{$4$} &
\multicolumn{1}{|c}{$3$} & $2$ & $4$ & \multicolumn{1}{|c}{$4$} & $4$ & $3$ &
$4$ & $3$ & $3$ & $4$ & \multicolumn{1}{c|}{$3$}\\\hline
\multicolumn{1}{|c}{$\{1,2,4\}$} & \multicolumn{1}{|c}{$4$} &
\multicolumn{1}{|c}{$4$} & $4$ & $4$ & \multicolumn{1}{|c}{$4$} & $4$ & $4$ &
$4$ & $4$ & $4$ & $4$ & \multicolumn{1}{c|}{$4$}\\\hline
\multicolumn{1}{|c}{$\{1,3,4\}$} & \multicolumn{1}{|c}{$4$} &
\multicolumn{1}{|c}{$4$} & $4$ & $4$ & \multicolumn{1}{|c}{$4$} & $4$ & $4$ &
$4$ & $4$ & $4$ & $4$ & \multicolumn{1}{c|}{$4$}\\\hline
\multicolumn{1}{|c}{$(4,2,3)$} & \multicolumn{1}{|c}{$4$} &
\multicolumn{1}{|c}{$4$} & $2$ & $3$ & \multicolumn{1}{|c}{$3$} & $3$ & $4$ &
$3$ & $4$ & $4$ & $3$ & \multicolumn{1}{c|}{$4$}\\
\multicolumn{1}{|c}{$(2,3,4)$} & \multicolumn{1}{|c}{$4$} &
\multicolumn{1}{|c}{$2$} & $3$ & $4$ & \multicolumn{1}{|c}{$3$} & $4$ & $3$ &
$3$ & $4$ & $3$ & $4$ & \multicolumn{1}{c|}{$4$}\\
\multicolumn{1}{|c}{$(3,4,2)$} & \multicolumn{1}{|c}{$4$} &
\multicolumn{1}{|c}{$3$} & $4$ & $2$ & \multicolumn{1}{|c}{$3$} & $3$ & $3$ &
$4$ & $4$ & $4$ & $4$ & \multicolumn{1}{c|}{$3$}\\\hline
\multicolumn{1}{|c}{$(2,4,3)$} & \multicolumn{1}{|c}{$4$} &
\multicolumn{1}{|c}{$2$} & $4$ & $3$ & \multicolumn{1}{|c}{$4$} & $3$ & $4$ &
$4$ & $3$ & $4$ & $3$ & \multicolumn{1}{c|}{$3$}\\
\multicolumn{1}{|c}{$(4,3,2)$} & \multicolumn{1}{|c}{$4$} &
\multicolumn{1}{|c}{$4$} & $3$ & $2$ & \multicolumn{1}{|c}{$4$} & $4$ & $4$ &
$3$ & $3$ & $3$ & $3$ & \multicolumn{1}{c|}{$4$}\\
\multicolumn{1}{|c}{$(3,2,4)$} & \multicolumn{1}{|c}{$4$} &
\multicolumn{1}{|c}{$3$} & $2$ & $4$ & \multicolumn{1}{|c}{$4$} & $4$ & $3$ &
$4$ & $3$ & $3$ & $4$ & \multicolumn{1}{c|}{$3$}\\\hline
&  &  &  &  &  &  &  &  &  &  &  &
\end{tabular}
\caption{{}{}Nonconservative majority operations generating a minimal clone on
the four-element set $\left\{  1,2,3,4\right\}  $}%
\label{table 4maj}%
\end{table}

Suppose that $f$ is a majority operation on $A$ generating a minimal clone. If
$f$ is \emph{conservative}, i.e., it preserves every subset of $A$, then the
restriction of $f$ to any three-element subset has to be
isomorphic\footnote{By a slight abuse of terminology, we say that operations
$f$ and $g$ defined on sets $A$ and $B$, respectively, are isomorphic, if the
algebras $\left(  A;f\right)  $ and $\left(  B;g\right)  $ are isomorphic.} to
one of the 12 functions in Table~\ref{table 3maj}, and $f$ is uniquely
determined by these restrictions. However, the converse is not true: given a
conservative majority operation $f$ whose restriction to every three-element
subset is isomorphic to one of these 12 functions, it is not guaranteed that
$\left[  f\right]  $ is a minimal clone. The appropriate necessary and
sufficient condition for the minimality was given by B.~Cs\'{a}k\'{a}ny in
\cite{Cscons}. It turns out that a minimal clone generated by a conservative
majority operation contains either 1, 3, 8 or 24 majority operations. If both
$m_{2}$ and $m_{3}$ appear among the restrictions of $f$, then all possible
pairs of majority operations from $\left[  m_{2}\right]  \times\left[
m_{3}\right]  $ appear as restrictions of compositions of $f$; this yields 24
majority functions.

The investigation of minimal clones on the four-element set carried out in
\cite{W4maj} did not give any new examples: if $f$ is a majority operation on
a four-element set generating a minimal clone, then $f$ is either
conservative, or it is isomorphic to one of the 12 functions shown in
Table~\ref{table 4maj}. (The middle two rows mean that the value of the
functions on $\left(  a,b,c\right)  $ is $4$ whenever $\left\{  a,b,c\right\}
=\left\{  1,2,4\right\}  $ or $\left\{  a,b,c\right\}  =\left\{
1,3,4\right\}  $.) Restricting these functions to $\left\{  2,3,4\right\}  $
we get (isomorphic copies of) the 12 functions of Table~\ref{table 3maj}, and,
in fact, this restriction is a clone isomorphism.

In all of the above examples, the clone contains 1, 3, 8 or 24 majority
operations, and actually the ternary part of the clone, as a Menger algebra,
is determined up to isomorphism by its size (see \cite{Wafew} for details).
This gives rise to the question whether this is always the case. Some modest
steps have been taken in \cite{Wafew} to give an affirmative answer to this
question. However, it turns out that the answer is negative:\ we will prove
the following theorem in Section~\ref{sect 26}.

\begin{theorem}
\label{thm 26}There exists a minimal clone with 26 majority operations.
\end{theorem}

The other main result of this paper is that the same \textquotedblleft
trick\textquotedblright\ that makes it possible to construct a minimal clone
with $24$ majority operations using $m_{2}$ and $m_{3}$ works with any
majority operation $f$ in place of $m_{3}$, provided that $f$ is cyclically
symmetric, i.e., $f$ satisfies the identity $f\left(  x_{1},x_{2}%
,x_{3}\right) \approx f\left(  x_{2},x_{3},x_{1}\right)  $.

\begin{theorem}
\label{thm 78}If there is a minimal clone with $n$ majority operations one of
which is cyclically symmetric, then there is a minimal clone with $3n$
majority operations.
\end{theorem}

Since the clone containing 26 majority operations that we present in
Section~\ref{sect 26} is generated by a cyclically symmetric majority
operation, the above theorem implies that there is a minimal clone with 78
majority operations. In sum, what we know about the number of majority
operations in a minimal clone is that it can be 1, 3, 8, 24, 26 or 78, but it
cannot be 2 or 4 (cf.~\cite{Wafew}). We do not know if there are infinitely
many such numbers, and we do not even know whether every minimal clone
contains only finitely many majority operations.

\section{Proof of Theorem~1\label{sect 26}}

Theorem~\ref{thm 26} is a result of a computer search: we checked for each
non-conservative cyclically symmetric majority operation $f$ on a five-element
set whether $\left[  f\right]  $ is a minimal clone or not. We considered only
cyclically symmetric functions, because the number of all majority operations
is so huge, that the problem seems to be inaccessible. Even for the cyclically
symmetric case, the task took several weeks on several computers. The outcome
is that except for one function (up to isomorphism and up to permutation of
variables) the minimal clones contain 1 or 8 majority operations, and the
structure of the ternary part is the same as that of $\left[  m_{1}\right]  $
or $\left[  m_{3}\right]  $ (see Table~\ref{table 3maj}). The exceptional
function is the function $f_{1}$ in Table~\ref{table 26}; it generates a
minimal clone with 26 majority operations. This was first proven by computer,
but it is possible to verify it by human reasoning as well (to be presented in
this section). However, we do not have a \textquotedblleft
human\textquotedblright\ proof for the fact that this is the only cyclically
symmetric majority operation on the five-element set that yields a new kind of
minimal clone.\renewcommand{\arraystretch}{1.3}

\begin{table}[ptb]
\centering%
\begin{tabular}
[c]%
{cccp{0.45cm}p{0.45cm}p{0.45cm}p{0.45cm}p{0.45cm}p{0.45cm}p{0.45cm}p{0.45cm}}%
\cline{4-11}
&  &  & \multicolumn{1}{|p{0.45cm}}{$f_{1}$} & $f_{2}$ & $g_{1}^{u,v}$ &
$g_{2}^{u,v}$ & $g_{3}^{u,v}$ & $g_{4}^{u,v}$ & $g_{5}^{u,v}$ &
\multicolumn{1}{p{0.45cm}|}{$g_{6}^{u,v}$}\\\cline{3-11}\cline{3-11}
&  & \multicolumn{1}{|c|}{$\left\{  0,1,\overline{1}\right\}  $} &
\multicolumn{1}{|p{0.45cm}}{$1$} & $1$ & $1$ & $1$ & $1$ & $1$ & $1$ &
\multicolumn{1}{p{0.45cm}|}{$1$}\\\cline{3-11}
&  & \multicolumn{1}{|c|}{$\left\{  2,1,\overline{1}\right\}  $} &
\multicolumn{1}{|p{0.45cm}}{$1$} & $1$ & $1$ & $1$ & $1$ & $1$ & $1$ &
\multicolumn{1}{p{0.45cm}|}{$1$}\\\cline{3-11}
&  & \multicolumn{1}{|c|}{$\left\{  0,2,\overline{2}\right\}  $} &
\multicolumn{1}{|p{0.45cm}}{$2$} & $2$ & $2$ & $2$ & $2$ & $2$ & $2$ &
\multicolumn{1}{p{0.45cm}|}{$2$}\\\cline{3-11}
&  & \multicolumn{1}{|c|}{$\left\{  1,2,\overline{2}\right\}  $} &
\multicolumn{1}{|p{0.45cm}}{$2$} & $2$ & $2$ & $2$ & $2$ & $2$ & $2$ &
\multicolumn{1}{p{0.45cm}|}{$2$}\\\cline{3-11}
&  & \multicolumn{1}{|c|}{$\left\{  \overline{1},2,\overline{2}\right\}  $} &
\multicolumn{1}{|p{0.45cm}}{$2$} & $2$ & $2$ & $2$ & $2$ & $2$ & $2$ &
\multicolumn{1}{p{0.45cm}|}{$2$}\\\hline
\multicolumn{1}{|c}{$\left(  0,\overline{1},2\right)  $} &
\multicolumn{1}{|c}{$\left(  0,1,\overline{2}\right)  $} &
\multicolumn{1}{|c|}{$\left(  0,1,2\right)  $} &
\multicolumn{1}{|p{0.45cm}}{$1$} & $2$ & $1$ & $2$ & $2$ & $2$ & $1$ &
\multicolumn{1}{p{0.45cm}|}{$1$}\\
\multicolumn{1}{|c}{$\left(  \overline{1},2,0\right)  $} &
\multicolumn{1}{|c}{$\left(  1,\overline{2},0\right)  $} &
\multicolumn{1}{|c|}{$\left(  1,2,0\right)  $} &
\multicolumn{1}{|p{0.45cm}}{$1$} & $2$ & $2$ & $1$ & $2$ & $1$ & $2$ &
\multicolumn{1}{p{0.45cm}|}{$1$}\\
\multicolumn{1}{|c}{$\left(  2,0,\overline{1}\right)  $} &
\multicolumn{1}{|c}{$\left(  \overline{2},0,1\right)  $} &
\multicolumn{1}{|c|}{$\left(  2,0,1\right)  $} &
\multicolumn{1}{|p{0.45cm}}{$1$} & $2$ & $2$ & $2$ & $1$ & $1$ & $1$ &
\multicolumn{1}{p{0.45cm}|}{$2$}\\\hline
\multicolumn{1}{|c}{$\left(  2,\overline{1},0\right)  $} &
\multicolumn{1}{|c}{$\left(  \overline{2},1,0\right)  $} &
\multicolumn{1}{|c|}{$\left(  2,1,0\right)  $} &
\multicolumn{1}{|p{0.45cm}}{$2$} & $1$ & $1$ & $2$ & $1$ & $2$ & $1$ &
\multicolumn{1}{p{0.45cm}|}{$2$}\\
\multicolumn{1}{|c}{$\left(  \overline{1},0,2\right)  $} &
\multicolumn{1}{|c}{$\left(  1,0,\overline{2}\right)  $} &
\multicolumn{1}{|c|}{$\left(  1,0,2\right)  $} &
\multicolumn{1}{|p{0.45cm}}{$2$} & $1$ & $1$ & $1$ & $2$ & $2$ & $2$ &
\multicolumn{1}{p{0.45cm}|}{$1$}\\
\multicolumn{1}{|c}{$\left(  0,2,\overline{1}\right)  $} &
\multicolumn{1}{|c}{$\left(  0,\overline{2},1\right)  $} &
\multicolumn{1}{|c|}{$\left(  0,2,1\right)  $} &
\multicolumn{1}{|p{0.45cm}}{$2$} & $1$ & $2$ & $1$ & $1$ & $1$ & $2$ &
\multicolumn{1}{p{0.45cm}|}{$2$}\\\hline
&  & \multicolumn{1}{|c|}{$\left(  0,\overline{1},\overline{2}\right)  $} &
\multicolumn{1}{|p{0.45cm}}{$\mathbf{1}$} & $\mathbf{2}$ & $\mathbf{u}$ &
$\mathbf{2}$ & $\mathbf{2}$ & $\mathbf{2}$ & $\mathbf{v}$ &
\multicolumn{1}{p{0.45cm}|}{$\mathbf{1}$}\\
&  & \multicolumn{1}{|c|}{$\left(  \overline{1},\overline{2},0\right)  $} &
\multicolumn{1}{|p{0.45cm}}{$\mathbf{1}$} & $\mathbf{2}$ & $\mathbf{2}$ &
$\mathbf{u}$ & $\mathbf{2}$ & $\mathbf{1}$ & $\mathbf{2}$ &
\multicolumn{1}{p{0.45cm}|}{$\mathbf{v}$}\\
&  & \multicolumn{1}{|c|}{$\left(  \overline{2},0,\overline{1}\right)  $} &
\multicolumn{1}{|p{0.45cm}}{$\mathbf{1}$} & $\mathbf{2}$ & $\mathbf{2}$ &
$\mathbf{2}$ & $\mathbf{u}$ & $\mathbf{v}$ & $\mathbf{1}$ &
\multicolumn{1}{p{0.45cm}|}{$\mathbf{2}$}\\\cline{3-11}
&  & \multicolumn{1}{|c|}{$\left(  \overline{2},\overline{1},0\right)  $} &
\multicolumn{1}{|p{0.45cm}}{$\mathbf{2}$} & $\mathbf{1}$ & $\mathbf{v}$ &
$\mathbf{2}$ & $\mathbf{1}$ & $\mathbf{2}$ & $\mathbf{u}$ &
\multicolumn{1}{p{0.45cm}|}{$\mathbf{2}$}\\
&  & \multicolumn{1}{|c|}{$\left(  \overline{1},0,\overline{2}\right)  $} &
\multicolumn{1}{|p{0.45cm}}{$\mathbf{2}$} & $\mathbf{1}$ & $\mathbf{1}$ &
$\mathbf{v}$ & $\mathbf{2}$ & $\mathbf{2}$ & $\mathbf{2}$ &
\multicolumn{1}{p{0.45cm}|}{$\mathbf{u}$}\\
&  & \multicolumn{1}{|c|}{$\left(  0,\overline{2},\overline{1}\right)  $} &
\multicolumn{1}{|p{0.45cm}}{$\mathbf{2}$} & $\mathbf{1}$ & $\mathbf{2}$ &
$\mathbf{1}$ & $\mathbf{v}$ & $\mathbf{u}$ & $\mathbf{2}$ &
\multicolumn{1}{p{0.45cm}|}{$\mathbf{2}$}\\\cline{3-11}
&  & \multicolumn{1}{|c|}{$\left(  \overline{2},1,\overline{1}\right)  $} &
\multicolumn{1}{|p{0.45cm}}{$1$} & $\overline{1}$ & $1$ & $\overline{1}$ &
$\overline{1}$ & $\overline{1}$ & $1$ & \multicolumn{1}{p{0.45cm}|}{$1$}\\
&  & \multicolumn{1}{|c|}{$\left(  1,\overline{1},\overline{2}\right)  $} &
\multicolumn{1}{|p{0.45cm}}{$1$} & $\overline{1}$ & $\overline{1}$ & $1$ &
$\overline{1}$ & $1$ & $\overline{1}$ & \multicolumn{1}{p{0.45cm}|}{$1$}\\
&  & \multicolumn{1}{|c|}{$\left(  \overline{1},\overline{2},1\right)  $} &
\multicolumn{1}{|p{0.45cm}}{$1$} & $\overline{1}$ & $\overline{1}$ &
$\overline{1}$ & $1$ & $1$ & $1$ & \multicolumn{1}{p{0.45cm}|}{$\overline{1}$%
}\\\cline{3-11}
&  & \multicolumn{1}{|c|}{$\left(  \overline{1},1,\overline{2}\right)  $} &
\multicolumn{1}{|p{0.45cm}}{$\overline{1}$} & $1$ & $1$ & $\overline{1}$ & $1
$ & $\overline{1}$ & $1$ & \multicolumn{1}{p{0.45cm}|}{$\overline{1}$}\\
&  & \multicolumn{1}{|c|}{$\left(  1,\overline{2},\overline{1}\right)  $} &
\multicolumn{1}{|p{0.45cm}}{$\overline{1}$} & $1$ & $1$ & $1$ & $\overline{1}
$ & $\overline{1}$ & $\overline{1}$ & \multicolumn{1}{p{0.45cm}|}{$1$}\\
&  & \multicolumn{1}{|c|}{$\left(  \overline{2},\overline{1},1\right)  $} &
\multicolumn{1}{|p{0.45cm}}{$\overline{1}$} & $1$ & $\overline{1}$ & $1$ & $1
$ & $1$ & $\overline{1}$ & \multicolumn{1}{p{0.45cm}|}{$\overline{1}$%
}\\\cline{3-11}
&  &  &  &  &  &  &  &  &  &
\end{tabular}
\caption{{}A minimal clone with 26 majority operations}%
\label{table 26}%
\end{table}

The base set for our functions will be $\left\{  0,1,2,\overline{1}%
,\overline{2}\right\}  $; this notation will help to emphasize certain
patterns in the functions. Table~\ref{table 26} shows the $26$ majority
operations in the clone under consideration. The functions are $f_{1}%
,f_{2},g_{1}^{u,v},\ldots,g_{6}^{u,v}$, where $u$ and $v$ can take the values
$1,\overline{1}$ independently of each other. Thus each column $g_{i}^{u,v}$
represents 4 functions, giving altogether $2+6\cdot4=26$ functions. The first
two functions are cyclically symmetric; we will use $f_{1}$ as a generator. We
only list the tuples where the majority rule does not apply, and in the first
five rows we make the same simplification as in Table~\ref{table 4maj}. For
example, the row $\left\{  0,1,\overline{1}\right\}  $ indicates that any one
of the 26 functions takes on the value $1$ on any triple $\left(
a,b,c\right)  $ such that $\left\{  a,b,c\right\}  =\left\{  0,1,\overline
{1}\right\}  $. Note also that for each of the functions the values on
$\left(  0,\overline{1},2\right)  ,\left(  0,1,\overline{2}\right)  $ and
$\left(  0,1,2\right)  $ coincide, and the same holds for permutations of
these triples.

The proof of Theorem~\ref{thm 26} consists of three lemmas. First we prove
that any majority operation generated by $f_{1}$ is one of the 26 functions
appearing in Table~\ref{table 26}, then we verify that these functions indeed
belong to $\left[  f_{1}\right]  $, and finally we prove that this clone is minimal.

\begin{lemma}
\label{lemma at most 26}The clone generated by $f_{1}$ contains at most $26$
majority operations.
\end{lemma}

\begin{proof}
Let $h$ be any majority function in $\left[  f_{1}\right]  $. We will examine
restrictions of $h$ to three- and four-element subsets in order to prove that
$h$ coincides with one of the 26 functions shown in Table~\ref{table 26}.

The three-element set $\left\{  0,1,\overline{1}\right\}  $ is preserved by
$f_{1}$, and the restriction of $f_{1}$ to this set is isomorphic to $m_{1}$.
There is only one majority function in $\left[  m_{1}\right]  $, therefore
$f_{1}$ and $h$ coincide on $\left\{  0,1,\overline{1}\right\}  $. The same
argument shows that the restriction of $h$ to any of the three-element sets
shown in the first five rows of the table is uniquely determined.

The four-element set $\left\{  \overline{1},0,2,1\right\}  $ is also preserved
by $f_{1}$, and the restriction to this set yields a function isomorphic to
$M_{3}$ (see Table~\ref{table 4maj}). This implies that there are 8
possibilities for $h$ on this four-element set, and $h|_{\left\{  \overline
{1},0,2,1\right\}  }$ is uniquely determined by $h|_{\left\{  0,1,2\right\}
}$. Similarly, there are 8 possibilities for $h$ on $\left\{  \overline
{2},0,1,2\right\}  $, and $h|_{\left\{  \overline{2},0,1,2\right\}  }$ is
uniquely determined by $h|_{\left\{  0,1,2\right\}  }$.

Let us also observe that $h$ preserves $\left\{  \overline{2},1,\overline
{1}\right\}  $, and its restriction to this set is isomorphic to $h|_{\left\{
0,1,2\right\}  }$, hence the latter determines $h|_{\left\{  \overline
{2},1,\overline{1}\right\}  }$.

We see that most values of $h$ are determined by $h|_{\left\{  0,1,2\right\}
} $, and the information we gathered about $h$ so far suffices to justify all
the entries in Table~\ref{table 26} except for the ones in boldface. Ignoring
these entries, i.e., the values on $\left\{  0,\overline{1},\overline
{2}\right\}  $, we have eight candidates for $h$, and the restrictions to
$\left\{  0,1,2\right\}  $ uniquely determine these (yet partial) functions.
Now we try to establish some relationships between the values on $\left\{
0,\overline{1},\overline{2}\right\}  $ and $\left\{  0,1,2\right\}  $. To this
end, we consider the smallest binary invariant relation$\footnote{In other
words, we consider the subalgebra generated by $\left\{  \left(  0,0\right)
,\left(  1,\overline{1}\right)  ,\left(  2,\overline{2}\right)  \right\}  $ in
the direct square of the algebra $\left(  \left\{  0,1,\overline
{1},2,\overline{2}\right\}  ;f_{1}\right)  $.}$ $\vartheta$ of $f_1$ relating
$0$ to $0$, $1$ to $\overline{1}$ and $2$ to $\overline{2}$:%
\[
\vartheta=\left\{
\begin{pmatrix}
0\\
0
\end{pmatrix}
,%
\begin{pmatrix}
1\\
1
\end{pmatrix}
,%
\begin{pmatrix}
1\\
\overline{1}%
\end{pmatrix}
,%
\begin{pmatrix}
2\\
2
\end{pmatrix}
,%
\begin{pmatrix}
2\\
\overline{2}%
\end{pmatrix}
\right\}  .
\]
Similarly, let $\varrho$ be the smallest invariant relation of $f_1$ that
relates $2$ to $0$, $1$ to $\overline{1}$ and $0$ to $\overline{2}$:%
\[
\varrho=\left\{
\begin{pmatrix}
0\\
\overline{2}%
\end{pmatrix}
,%
\begin{pmatrix}
1\\
1
\end{pmatrix}
,%
\begin{pmatrix}
1\\
\overline{1}%
\end{pmatrix}
,%
\begin{pmatrix}
1\\
2
\end{pmatrix}
,%
\begin{pmatrix}
2\\
0
\end{pmatrix}
,%
\begin{pmatrix}
2\\
1
\end{pmatrix}
,%
\begin{pmatrix}
2\\
2
\end{pmatrix}
\right\}  .
\]
We show below how the first one of the rows containing the boldface entries
can be filled out with the help of these relations; the other five rows can be
treated similarly.

Since $h$ belongs to $\left[  f_{1}\right]  $, it must preserve $\vartheta$
and $\varrho$, thus $h\left(  0,1,2\right)  \vartheta h\left(  0,\overline
{1},\overline{2}\right)  $ and $h\left(  2,1,0\right)  \varrho h\left(
0,\overline{1},\overline{2}\right)  $. We already know that $h\left(
0,1,2\right)  ,h\left(  2,1,0\right)  \in\left\{  1,2\right\}  $, so we have
the following four cases:%
\[
h\left(  0,\overline{1},\overline{2}\right)  =\left\{  \!\!\!%
\begin{array}
[c]{rl}%
2,~ & \text{if }h\left(  0,1,2\right)  =2,\,h\left(  2,1,0\right)  =2;\\
2,~ & \text{if }h\left(  0,1,2\right)  =2,\,h\left(  2,1,0\right)  =1;\\
1,~ & \text{if }h\left(  0,1,2\right)  =1,\,h\left(  2,1,0\right)  =2;\\
1\text{ or }\overline{1},~ & \text{if }h\left(  0,1,2\right)  =1,\,h\left(
2,1,0\right)  =1.
\end{array}
\right.
\]
We see that the value of $h\left(  0,\overline{1},\overline{2}\right)  $ is
uniquely determined except for two of the eight possibilities (denoted by
$\mathbf{u}$ and $\mathbf{v}$ in the table), where the value can be either $1$
or $\overline{1}$.
\end{proof}

\begin{lemma}
\label{lemma at least 26}The clone generated by $f_{1}$ contains at least $26$
majority operations.
\end{lemma}

\begin{proof}
We claim that the 26 functions shown in Table~\ref{table 26} belong to the
clone generated by $f_{1}$. The function $f_{2}$ can be obtained from $f_{1}$
by permuting variables, and, similarly, $g_{i}^{u,v}$ can be obtained from
$g_{1}^{u,v}$ for $i=2,\ldots,6$. Thus it suffices to show that $g_{1}%
^{u,v}\in\left[  f_{1}\right]  $ for all $u,v\in\left\{  1,\overline
{1}\right\}  $. This can be done by presenting a suitable composition for each
of these four functions:%
\begin{align*}
g_{1}^{1,1}\left(  x_{1},x_{2},x_{3}\right)   &  =f_{1}\left(  x_{2}%
,f_{1}\left(  x_{2},x_{1},x_{3}\right)  ,f_{1}\left(  x_{1},x_{2}%
,x_{3}\right)  \right)  ,\\
g_{1}^{1,\overline{1}}\left(  x_{1},x_{2},x_{3}\right)   &  =f_{1}\left(
x_{1},x_{2},f_{1}\left(  x_{2},x_{1},x_{3}\right)  \right)  ,\\
g_{1}^{\overline{1},1}\left(  x_{1},x_{2},x_{3}\right)   &  =f_{1}\left(
x_{3},x_{2},f_{1}\left(  x_{2},f_{1}\left(  x_{2},x_{1},x_{3}\right)
,f_{1}\left(  x_{1},x_{2},x_{3}\right)  \right)  \right)  ,\\
g_{1}^{\overline{1},\overline{1}}\left(  x_{1},x_{2},x_{3}\right)   &
=f_{1}\left(  x_{3},x_{2},f_{1}\left(  x_{1},x_{2},f_{1}\left(  x_{2}%
,x_{1},x_{3}\right)  \right)  \right)  .\qedhere
\end{align*}

\end{proof}

\begin{lemma}
\label{lemma minimal 26}The clone generated by $f_{1}$ is minimal.
\end{lemma}

\begin{proof}
We need to verify that (\ref{minimality criterion maj}) holds for $f_{1}$,
i.e., we have to prove that each of the 26 majority functions in $\left[
f_{1}\right]  $ generates $f_{1}$. Up to permutation of variables we have only
the five functions $f_{1},g_{1}^{1,1},g_{1}^{1,\overline{1}},g_{1}%
^{\overline{1},1},g_{1}^{\overline{1},\overline{1}}$. For the first one our
task is void, for the remaining four ones the same composition works: for all
$u,v\in\left\{  1,\overline{1}\right\}  $ we have%
\[
f_{1}\left(  x_{1},x_{2},x_{3}\right)  =g_{1}^{u,v}\left(  g_{1}^{u,v}\left(
x_{2},x_{1},x_{3}\right)  ,g_{1}^{u,v}\left(  x_{1},x_{3},x_{2}\right)
,g_{1}^{u,v}\left(  x_{3},x_{2},x_{1}\right)  \right)  .\qedhere
\]

\end{proof}

\section{Proof of Theorem~2\label{sect 78}}

Let $f$ be a cyclically symmetric majority operation on a set $A$ that
generates a minimal clone containing $n$ majority operations. We add a new
element to the base set: $A^{\ast}\mathrel{\mathop:}=A\dot{\cup}\left\{
\ast\right\}  $, and we construct a majority function $f^{\ast}$ on $A^{\ast}$
as follows. If $a_{1},a_{2},a_{3}\in A^{\ast}$ are not pairwise different,
then we define $f^{\ast}\left(  a_{1},a_{2},a_{3}\right)  $ by the majority
rule, otherwise let%
\[
f^{\ast}\left(  a_{1},a_{2},a_{3}\right)  =\left\{  \!\!\!%
\begin{array}
[c]{rl}%
f\left(  a_{1},a_{2},a_{3}\right)  ,~ & \text{if }\left\{  a_{1},a_{2}%
,a_{3}\right\}  \subseteq A;\\
a_{1},~ & \text{if }\ast\in\left\{  a_{1},a_{2},a_{3}\right\}  .
\end{array}
\right.
\]
Let us observe that $f^{\ast}|_{A}$ coincides with $f$, and $f^{\ast}|_{B}$ is
isomorphic to $m_{2}$ for any three-element set $B\subseteq A^{\ast}$ that is
not a subset of $A$. We claim that $f^{\ast}$ generates a minimal clone with
$3n$ majority operations. Just as in the previous section, we divide the task
into three lemmas. First we prove that $3n$ is an upper bound for the number
of majority operations in $\left[  f^{\ast}\right]  $, then we show that this
bound is sharp, and, finally, we verify that the clone is minimal.

\begin{lemma}
\label{lemma at most 3n}The clone generated by $f^{\ast}$ contains at most
$3n$ majority operations.
\end{lemma}

\begin{proof}
Let $h^{\ast}$ be any majority function in $\left[  f^{\ast}\right]  $. Since
$f^{\ast}$ preserves $A$ and all three-element subsets $B$ with $\ast\in B$,
the function $h^{\ast}$ must preserve these sets as well. Clearly, $h^{\ast}$
is uniquely determined by its restrictions to all previously mentioned sets.
The restriction of $h^{\ast}$ to $A$ belongs to $\left[  f\right]  $, hence
there are $n$ possibilities for~$h^{\ast}|_{A}$. If $B$ is a three-element set
containing $\ast$, then $h^{\ast}|_{B}$ is isomorphic to one of the three
majority functions in $\left[  m_{2}\right]  $. Moreover, if $B_{1}$ and
$B_{2}$ are two such subsets, then $h^{\ast}|_{B_{1}}$ and $h^{\ast}|_{B_{2}}$
are isomorphic, since the same holds for $f^{\ast}$. This means that we cannot
choose $h^{\ast}|_{B_{1}}$ and $h^{\ast}|_{B_{2}}$ independently: if one of
them is given, the other one is uniquely determined. Thus $h^{\ast}$ is
determined by $h^{\ast}|_{A}$ and $h^{\ast}|_{B}$ for \emph{one} three-element
set $B$ with $\ast\in B$, hence there are (at most) $3n$ possibilities for the function
$h^{\ast}$.
\end{proof}

For the rest of the paper it will be convenient to introduce some notation.
Let us rename the function $m_{2}$ to $d_{1}$, and let us denote the other two
majority operations in its clone by $d_{2}$ and $d_{3}$ (see
Table~\ref{table 3maj}). The motivation for the notation is that $d_{i}$
coincides with the $i$th projection whenever its arguments are pairwise
different:%
\[
d_{i}\left(  a_{1},a_{2},a_{3}\right)  =a_{i}\text{ for }\left\{  a_{1}%
,a_{2},a_{3}\right\}  =\left\{  1,2,3\right\}  \text{.}%
\]
Observe also that $d_{3}$ is the dual discriminator function on $\left\{
1,2,3\right\}  $.

We have seen in the proof of the above lemma, that for any majority operation
$h^{\ast}\in\left[  f^{\ast}\right]  $, the restriction of $h^{\ast}$ to $A$
is a majority function $h\in\left[  f\right]  $, and the restriction of
$h^{\ast}$ to any three-element subset $B$ that contains $\ast$ is isomorphic
to $d_{i}$ for some $i\in\left\{  1,2,3\right\}  $ (where $i$ does not depend
on $B$). Since $h^{\ast}$ is determined by $h$ and $d_{i}$, we will use the
notation $h^{\ast}=h\ast d_{i}$. For example, we have $f^{\ast}=f\ast d_{1}$.
Equivalently, $h\ast d_{i}$ is the majority operation on $A^{\ast}$ defined
for pairwise different $a_{1},a_{2},a_{3}\in A^{\ast}$ by%
\[
\left(  h\ast d_{i}\right)  \left(  a_{1},a_{2},a_{3}\right)  =\left\{  \!\!\!%
\begin{array}
[c]{rl}%
h\left(  a_{1},a_{2},a_{3}\right)  ,~ & \text{if }\left\{  a_{1},a_{2}%
,a_{3}\right\}  \subseteq A;\\
a_{i},~ & \text{if }\ast\in\left\{  a_{1},a_{2},a_{3}\right\}  .
\end{array}
\right.
\]

In the following claim, which is the key for proving that $\left[  f^{\ast
}\right]  $ contains at least $3n$ majority operations, we will consider terms
involving a ternary operation symbol $d$ and the variables $x_{1},x_{2},x_{3}
$. We say that a term $s$ \emph{is obtained from the term }$t$\emph{\ using
cyclic shifts}, iff $s$ arises from $t$ by a finite number of replacements of
some subterm $d\left(  t_{1},t_{2},t_{3}\right)  $ of $t$ by $d\left(
t_{2},t_{3},t_{1}\right)  $ or $d\left(  t_{3},t_{1},t_{2}\right)  $. More
formally, the \emph{set }$\operatorname*{CS}$\emph{$(t)$ of all terms that can
be obtained from $t$ using cyclic shifts} is defined inductively: If $t$ is a
variable, then $\operatorname*{CS}(t)\mathrel{\mathop:}=\left\{  t\right\}  $.
Otherwise, if $t=d\left(  t_{1},t_{2},t_{3}\right)  $ where the sets
$\operatorname*{CS}(t_{i})$, $i=1,2,3$, are already defined, then
\[
\operatorname*{CS}(t)\mathrel{\mathop:}=\left\{  \left.  d\left(  s_{1}%
,s_{2},s_{3}\right)  ,d\left(  s_{2},s_{3},s_{1}\right)  ,d\left(  s_{3}%
,s_{1},s_{2}\right)  \ \right\vert \ s_{i}\in\operatorname*{CS}(t_{i})\text{
for }i=1,2,3\right\}  .
\]
A straightforward induction argument shows that if one evaluates such terms
over an algebra with a cyclically symmetric basic operation, then every term
has the same term function as all its cyclically shifted descendants.

We are going to evaluate our terms over the algebra $\mathbb{B}=\left(
\left\{  1,2,3\right\}  ;d_{1}\right)  $. The term function corresponding to
the term $t$ is denoted by $t^{\mathbb{B}}$. With this notation we have
$d_{1}=\left(  d\left(  x_{1},x_{2},x_{3}\right)  \right)  ^{\mathbb{B}}$, and
$x_{i}^{\mathbb{B}}$ is the $i$th ternary projection on $\left\{
1,2,3\right\}  $. (Let us note that we did not distinguish variables and
projections until now.)

\begin{claim}
\label{claim cyclic shifts}If $t$ is a ternary term over the algebra
$\mathbb{B}=\left(  \left\{  1,2,3\right\}  ;d_{1}\right)  $ and
$t^{\mathbb{B}}$ is not a projection, then there exist terms $s_{1}%
,s_{2},s_{3}\in\operatorname*{CS}\left(  t\right)  $ such that%
\[
s_{1}^{\mathbb{B}}=d_{1},~s_{2}^{\mathbb{B}}=d_{2}~\text{and~}s_{3}%
^{\mathbb{B}}=d_{3}.
\]

\end{claim}

\begin{proof}
We prove the claim by term induction. Since $t$ does not evaluate to a
projection, the term function corresponding to $t$ is a majority function. It
is not hard to show (e.g., by another term induction) that after performing
arbitrary cyclic shifts, we still get a majority operation, i.e., one of the
functions $d_{1},d_{2},d_{3}$. We need to show that we can actually get all
three of these functions.

Since $t$ cannot be a single variable, the initial step of the induction is
the case when it contains only one operation symbol $d$. Then $t$ is of the
form $d\left(  x_{i_{1}},x_{i_{2}},x_{i_{3}}\right)  $. The indices
$i_{1},i_{2},i_{3}$ must be pairwise different (otherwise $t$ would evaluate
to a projection), thus $t^{\mathbb{B}}=d_{i_{1}}$. Using cyclic shifts we get
$\left(  d\left(  x_{i_{2}},x_{i_{3}},x_{i_{1}}\right)  \right)  ^{\mathbb{B}%
}=d_{i_{2}}$ and $\left(  d\left(  x_{i_{3}},x_{i_{1}},x_{i_{2}}\right)
\right)  ^{\mathbb{B}}=d_{i_{3}}$. Since $\left\{  i_{1},i_{2},i_{3}\right\}
=\left\{  1,2,3\right\}  $, we have obtained all three of $d_{1},d_{2},d_{3}$
(in some order).

For the inductive step let us write $t$ in the form $t=d\left(  t_{1}%
,t_{2},t_{3}\right)  $. For each $k=1,2,3$ we have two possibilities for the
subterm $t_{k}$: either $t_{k}^{\mathbb{B}}=x_{i_{k}}^{\mathbb{B}}$ or
$t_{k}^{\mathbb{B}}=d_{i_{k}}$ for some $i_{k}\in\left\{  1,2,3\right\}  $. By
induction, in the latter case we can obtain any one of $d_{1},d_{2},d_{3}$ by
performing cyclic shifts on $t_{k}$.

Suppose that some of $i_{1},i_{2},i_{3}$ coincide, say $i_{1}=i_{2}\neq i_{3}
$. If $t_{1}^{\mathbb{B}}=t_{2}^{\mathbb{B}}=x_{i_{1}}^{\mathbb{B}}$, then
$t^{\mathbb{B}}=\left(  d\left(  x_{i_{1}},x_{i_{1}},t_{3}\right)  \right)
^{\mathbb{B}}=x_{i_{1}}^{\mathbb{B}}$ by the majority rule, contradicting our
assumption that $t$ does not evaluate to a projection. Thus at least one of
$t_{1}^{\mathbb{B}},t_{2}^{\mathbb{B}}$ equals $d_{i_{1}}$, say $t_{1}%
^{\mathbb{B}}=d_{i_{1}}$. Let us choose $j_{1}\in\left\{  1,2,3\right\}  $ to
be different from $i_{2}$ and $i_{3}$, and construct a term $t_{1}^{\prime}$
using cyclic shifts in $t_{1}$ (by induction), such that $t_{1}^{\prime
\mathbb{B}}=d_{j_{1}}$. For compatibility, let us put $j_{2}=i_{2},j_{3}%
=i_{3}$ and $t_{2}^{\prime}=t_{2},t_{3}^{\prime}=t_{3}$. The key property of
the new terms $t_{1}^{\prime},t_{2}^{\prime},t_{3}^{\prime}$ is the following:%
\begin{equation}
t_{k}^{\prime}\in\operatorname*{CS}\left(  t_{k}\right)  ,~t_{k}%
^{\prime\mathbb{B}}=x_{j_{k}}^{\mathbb{B}}\text{ or }t_{k}^{\prime\mathbb{B}%
}=d_{j_{k}}\text{, and }\left\{  j_{1},j_{2},j_{3}\right\}  =\left\{
1,2,3\right\}  \text{.} \label{key property}%
\end{equation}
If $i_{1}=i_{2}=i_{3}$, then at most one of $t_{1},t_{2},t_{3}$ can evaluate
to a projection, and applying cyclic shifts on the other two terms we can
still achieve (\ref{key property}). If $i_{1},i_{2},i_{3}$ are pairwise
different, then (\ref{key property}) is trivially achieved by letting
$j_{k}=i_{k}$ and $t_{k}^{\prime}=t_{k}$ for $k=1,2,3$.

Now we are in a position to write up the desired terms $s_{1},s_{2},s_{3}$
(keeping in mind that $\left\{  j_{1},j_{2},j_{3}\right\}  =\left\{
1,2,3\right\}  $):%
\begin{align*}
s_{j_{1}}\mathrel{\mathop:}  &  =d\left(  t_{1}^{\prime},t_{2}^{\prime}%
,t_{3}^{\prime}\right)  ;\\
s_{j_{2}}\mathrel{\mathop:}  &  =d\left(  t_{2}^{\prime},t_{3}^{\prime}%
,t_{1}^{\prime}\right)  ;\\
s_{j_{3}}\mathrel{\mathop:}  &  =d\left(  t_{3}^{\prime},t_{1}^{\prime}%
,t_{2}^{\prime}\right)  .
\end{align*}

These terms can be obtained from $t$ by cyclic shifts, therefore, the
corresponding term functions are majority functions. Thus it suffices to
verify the equalities $s_{j_{1}}^{\mathbb{B}}=d_{j_{1}},s_{j_{2}}^{\mathbb{B}%
}=d_{j_{2}},s_{j_{3}}^{\mathbb{B}}=d_{j_{3}}$ for tuples $\left(  a_{1}%
,a_{2},a_{3}\right)  \in\left\{  1,2,3\right\}  ^{3}$ where $a_{1},a_{2}%
,a_{3}$ are pairwise different. For such a tuple we have $t_{k}^{\prime
\mathbb{B}}\left(  a_{1},a_{2},a_{3}\right)  =a_{j_{k}}$, regardless of
whether $t_{k}^{\prime\mathbb{B}}=x_{j_{k}}^{\mathbb{B}}$ or $t_{k}%
^{\prime\mathbb{B}}=d_{j_{k}}$, hence%
\begin{align*}
s_{j_{1}}^{\mathbb{B}}\left(  a_{1},a_{2},a_{3}\right)   &  =d_{1}\left(
a_{j_{1}},a_{j_{2}},a_{j_{3}}\right)  =a_{j_{1}}=d_{j_{1}}\left(  a_{1}%
,a_{2},a_{3}\right) \\
s_{j_{2}}^{\mathbb{B}}\left(  a_{1},a_{2},a_{3}\right)   &  =d_{1}\left(
a_{j_{2}},a_{j_{3}},a_{j_{1}}\right)  =a_{j_{2}}=d_{j_{2}}\left(  a_{1}%
,a_{2},a_{3}\right) \\
s_{j_{3}}^{\mathbb{B}}\left(  a_{1},a_{2},a_{3}\right)   &  =d_{1}\left(
a_{j_{3}},a_{j_{1}},a_{j_{2}}\right)  =a_{j_{3}}=d_{j_{3}}\left(  a_{1}%
,a_{2},a_{3}\right)  .\qedhere
\end{align*}

\end{proof}

Let us observe that we did not really use the fact that $\mathbb{B}$ has only
three elements: the claim is true for $\mathbb{B}=\left(  B;d_{1}\right)  $
for an arbitrary nonempty set $B$, where $d_{1}$ is the majority operation
defined for pairwise distinct $a_{1},a_{2},a_{3}\in B$ by the formula
$d_{1}\left(  a_{1},a_{2},a_{3}\right)  =a_{1}$.

\begin{lemma}
\label{lemma at least 3n}The clone generated by $f^{\ast}$ contains at least
$3n$ majority operations.
\end{lemma}

\begin{proof}
We will prove that for any majority function $h\in\left[  f\right]  $, the
three functions $h\ast d_{1}$, $h\ast d_{2}$, $h\ast d_{3}$ belong to $\left[
f^{\ast}\right]  $. Since $h\in\left[  f\right]  $, there is a composition of
$f$ that yields $h$. This composition can be described by a ternary term $t$
such that the corresponding term function over the algebra $\left(
A;f\right)  $ is $h$. Since $h$ is a majority operation, the term operation
$t^{\mathbb{B}}$ over the algebra $\mathbb{B}\mathrel{\mathop:}=(\left\{
1,2,3\right\}  ;d_{1})$ is not a projection, and
Claim~\ref{claim cyclic shifts} is applicable.

Let $s_{1},s_{2},s_{3}$ be the terms constructed from $t$ by cyclic shifts in
Claim~\ref{claim cyclic shifts}, and let $h_{1}^{\ast},h_{2}^{\ast}%
,h_{3}^{\ast}$ be the corresponding term functions over the algebra $\left(
A^{\ast};f^{\ast}\right)  $; these functions clearly belong to $\left[
f^{\ast}\right]  $. If $B$ is a three-element subset of $A^{\ast}$ that is not
contained in $A$, then $h_{i}^{\ast}|_{B}$ is isomorphic to $d_{i}$ as
$(B;f^{\ast}|_{B})\cong\mathbb{B}$. Since $f^{\ast}|_{A}=f$ is cyclically
symmetric, the cyclic shifts do not change the term functions on $A$: we have
$h_{i}^{\ast}|_{A}=h$. Thus we can conclude that $h_{i}^{\ast}=h\ast d_{i}$
for $i=1,2,3$.

Since there are $n$ choices for $h$, the clone generated by $f^{\ast}$
contains the $3n$ functions $h\ast d_{i}\left(  h\in\left[  f\right]
,i\in\left\{  1,2,3\right\}  \right)  $.
\end{proof}

\begin{remark}
The{\ previous two lemmas can be interpreted from the viewpoint of abstract
clones as follows. }For{\ a fixed subset $B\subseteq A^{\ast}$ with }$\ast\in
B${, the restriction mappings}%
\begin{align*}
|_{A}\colon{[f^{\ast}]}^{\left(  3\right)  }  &  \rightarrow{[f^{\ast}%
|_{A}]^{\left(  3\right)  }=}\left[  f\right]  ^{\left(  3\right)  }\\
|_{B}\colon{[f^{\ast}]}^{\left(  3\right)  }  &  \rightarrow{[f^{\ast}|_{B}%
]}^{\left(  3\right)  }{\cong}\left[  m_{2}\right]  ^{\left(  3\right)  }%
\end{align*}
{are homomorphisms of Menger algebras (since they are induced by clone
homomorphisms). The proof of Lemma~\ref{lemma at most 3n} shows that the
intersection of the kernels of these two homomorphisms is the equality
relation on }${[f^{\ast}]}^{\left(  3\right)  }$, and the proof of
{Lemma~\ref{lemma at least 3n} shows that these homomorphisms are surjective.
Thus the Menger algebra $[f^{\ast}]^{(3)}$, which is decisive for the
minimality of $[f^{\ast}]$, is a subdirect product of the Menger algebras $[f]^{(3)}$ and~
$[m_{2}]^{(3)}$.}
\end{remark}

\begin{lemma}
\label{lemma minimal 3n}The clone generated by $f^{\ast}$ is minimal.
\end{lemma}

\begin{proof}
According to (\ref{minimality criterion maj}), we need to prove that for any
majority operation $h^{\ast}$ in the clone generated by $f^{\ast}$, we have
$f^{\ast}\in\left[  h^{\ast}\right]  $. We know that $h^{\ast}$ is of the form
$h\ast d_{i}$, where $h\in\left[  f\right]  $ and $i\in\left\{  1,2,3\right\}
$. Since the clone generated by $f$ is minimal, there is a composition that
produces $f$ from $h$. Applying this composition for $h^{\ast}=h\ast d_{i}$,
we get a function of the form $f\ast d_{j}$. Taking into account that $f$ is
cyclically symmetric, a suitable cyclic permutation of variables of $f\ast
d_{j}$ yields $f^{\ast}=f\ast d_{1}$.
\end{proof}

\section*{Acknowledgements}

The authors would like to thank everyone that directly or indirectly
contributed to the obtainment of the results presented in this article. In
fact, these would not exist without many people that assisted in the computer
search for cyclically symmetric minimal majority functions on a five-element
set, and only few of which can be given credit here in a personal form.

In particular, the authors would like to mention G\'{e}za Makay for his
recommendations concerning the design of the computer program that first found
the minimal clone in Table~\ref{table 26}. Furthermore, they express their
sincere gratitude towards the staff of the Department of Mathematics at
Dresden University of Technology in general for the generous provision of
computing capacities, and to Mr.~Pierre Frison, who helped by running the
program on his private laptop.

The second named author acknowledges that the present project is supported by
the National Research Fund, Luxembourg, and cofunded under the Marie Curie
Actions of the European Commission (FP7-COFUND), and supported by the
Hungarian National Foundation for Scientific Research under grants no.~K60148
and K77409.

\end{document}